\documentstyle[12pt]{article}

\newtheorem{theorem}{Theorem}
\newtheorem{corollary}{Corollary}

\newtheorem{lemma}{Lemma}

\newcommand{\R}{\hbox{\rm\setbox1=\hbox{I}\copy1\kern-.45\wd1 R}}

\newcommand{\r}{\scriptsize{\hbox{\rm\setbox1=\hbox{I}\copy1\kern-.45\wd1 R}}}

\newcommand{\proof}{\medskip  \par \noindent {\bf Proof.} \ \ }

\newcommand{\bx}{{\bf X} }

\newcommand{\p}{{\cal P}}

\begin{document}
\bibliographystyle{plain}

\thispagestyle{empty}
\setcounter{page}{0}

\vspace {2cm}

{\large G. Morvai, S. Yakowitz, and L. Gy\"orfi: }

\vspace {1cm}

{\Large Nonparametric inference for ergodic, stationary time series.}

\vspace {1cm}

{\large   Ann. Statist.  24  (1996),  no. 1, 370--379.}

\vspace {2cm}

\begin{abstract}
The setting is a stationary, ergodic time series.
The challenge is to construct a sequence of functions, each based on
only  finite segments of the past, which
together provide a strongly consistent estimator for the
conditional probability of the next observation, given the infinite
past.  Ornstein gave such a construction for
the  case that the  values are from a finite set,
and recently Algoet extended the
scheme to  time series with coordinates in a
Polish space.

 The present study relates a different solution
to the challenge.  The algorithm is simple and its verification is
fairly transparent.  Some extensions to regression, pattern
recognition,
and on-line forecasting are mentioned.

\end{abstract}

\pagebreak

\section{Introduction}
\label{sec1}

In this section, we give brief overview of the situation
with respect to nonparametric inference under the most
lenient mixing conditions.
 Impetus for this line of study follows
Roussas (1969) and Rosenblatt (1970) who extended ideas in
the nonparametric regression literature for i.i.d. variables to
give  a theory adequate for showing, for example,
that for $\{X_i\}$ a real Markov sequence, under Doeblin-like
assumptions, the obvious kernel forecaster is an asymptotically
normal estimator of the conditional expectation
$E(X_0|X_{-1}=x)$. In the 1980's, there was an explosion of works
 which showed consistency in various senses for nonparametric
auto-regression and density estimators under more and more general mixing
assumptions (e.g., Castellana and Leadbetter (1986), Collomb (1985),
Gy\"orfi (1981), and Masry (1986)).
The monograph by Gy\"orfi {\it et
al.} (1989) gives supplemental information about nonparametric
estimation for dependent series.

Such striving for generality stems from  the
inconvenience of mixing conditions; satisfactory statistical
tests are not available.  Some recent developments have succeeded
in disposing of these conditions altogether.  In the
Markov case, aside from some smoothness assumptions,
it is enough that an invariant law exist to get
the usual pointwise asymptotic normality of kernel regression (Yakowitz
(1989)).  In case of Harris recurrence but  no invariant law,
one can still attain  a.s. pointwise convergence of a nearest-neighbor
regression algorithm in which the neighborhood is chosen in
advance
and observations continue until
a prescribed number of points fall into that neighborhood (Yakowitz (1993)).

Pushing beyond the Markov
hypothesis, by a histogram estimate (Gy\"orfi {\it et
al.} (1989)) or a  recursive-type estimator
(Gy\"orfi and Masry (1990)), one can infer the marginal
density of an ergodic stationary time series provided only
that there exist an absolutely continuous transition density.
Here the limit may have been attained;  it is now known
(Gy\"orfi {\it et al.} (1989) and Gy\"orfi and Lugosi (1992),
respectively) that without the conditional density assumption,
the histogram estimator and the
 kernel and recursive kernel estimates
for the marginal density are not generally consistent. 

The situation with respect to (auto-) regression is more inclusive
for ergodic, stationary sequences.
In a landmark paper, following developments by Ornstein (1978)
for the case   that the time series values
are from a finite set, for time series with values
in a Polish space, Algoet (1992, \S 5) has  provided a data-driven
distribution function construction $F_n(x|X_{-1},X_{-2},\dots)$ which a.s. converges in
distribution to
$$P(X_0\le x|X_{-1},X_{-2},\dots)=P(X_{0}\le x|\bx^-),$$
where ${\bf X^-} = {(X_{-1},X_{-2},\dots)}$.

 The goal of the present study is to relate a
simpler rule the consistency of which is easy to establish.
In concluding sections,
it is noted that as a result of these developments, one has
a consistent regression estimate in the bounded time-series case,
and implications to problems of pattern recognition and
 on-line forecasting are mentioned.  It is to be conceded that our algorithm,
as well as those of Algoet's and Ornstein's, can be expected to
require very large data segments for acceptable precision.

As a final general comment, we note that the assumption
of ergodicity may be relaxed somewhat.  Thus
in view of Sections 7.4 and 8.5 of Gray (1988), one sees that
a nonergodic stationary process has an ergodic decomposition.
With probability one, a realization of the time series falls into
an invariant event on which the process is ergodic and
stationary. Then one may apply  the developments of this study to
that event as though it were the process universe. Thus the analysis
here also remains valid for stationary nonergodic processes.
Our analysis is restricted to the case that the coordinates
of the time series are real, but it is evident that the proofs
extend directly to the vector-valued case.  In view of
Theorem 2.2 of Billingsley (1968, p. 14) it will be clear that
the formulas and derivations to follow also hold if the $X_i's$
are in a Polish space.

\section{Estimation of conditional distributions}
\label{sec2}

Let $\bx =\{X_n\}$ denote a real-valued doubly
infinite  stationary ergodic time series.
Let $$X^{-1}_{-j}=(X_{-j},X_{-j+1},\dots,X_{-1})$$
be notation for a data segment into the j-past, where j
may be infinite.
For a Borel set $C$ one wishes to infer  the conditional
probability
$$P(C|\bx ^-)=P(X_0\in C|X_{-\infty}^{-1}). $$

The algorithm to be
promoted here is iterative on an index $k=1,2,\dots$  For
each $k$, the data-driven estimate  of
$P(C|\bx ^-)$ requires only a segment of finite (but random)
length
of $\bx ^-$. One may proceed by simply repeating the estimation
process for $k$=1,2,\dots, until a given finite data record no longer
suffices for the demands of the algorithm. The goal of the study
will be to show that
a.s. convergence can be attained. That is, our estimation is
strongly consistent in the topology of weak convergence.

The estimation algorithm is now revealed in the simple context
of binary sequences,
 and afterwards, we show alterations necessary
for more general processes.

Define the sequences $\lambda _{k-1}$ and $\tau _k$ recursively ($k=1,2,\dots$).
Put $\lambda _0=1$ and let $\tau _k$ be the time between the 
occurrence of the pattern 
$$B(k)=(X_{-\lambda_{k-1}},\ldots,X_{-1})=
X_{-\lambda_{k-1}}^{-1}
$$
at time $-1$ and the last occurrence of the same pattern prior to time 
$-1$. More precisely, let
$$
{\tau}_k=
\min\{t>0 : X_{-\lambda_{k-1}-t}^{-1-t}=X_{-\lambda_{k-1}}^{-1}\}.
$$
Put
$$\lambda _k=\tau _k+\lambda_{k-1}.$$
The observed vector
$B(k)$ a.s. takes a value  having positive probability; thus  by ergodicity,
with probability $1$
the string $B(k)$ must appear infinitely often in the sequence
$X_{-\infty}^{-2}$.
One denotes the $k$th estimate of $P(C|\bx ^-)$ by $P_k(C),$
and defines it to be
\begin{equation}
\label{eq1}
P_k(C)={1\over k}\sum_{1\le j\le k} 1_{C}(X_{-\tau _j}).
\end{equation}
Here $1_{C}$ is the indicator function for $C$.

For the general case,
we use a sub-sigma-field structure motivated by Algoet (1992, Section
5.2), which is more general. Let
${\p}_k=\{ A_{k,i},\, i=1,2,\dots ,m_k\}$
 be a sequence of finite
partitions of the real line by (finite or infinite) right semi-closed 
intervals such that $\sigma ({\p}_k)$ is an increasing sequence 
of finite $\sigma$-algebras 
that asymptotically generate the Borel $\sigma$-field.
Let $G_k$ denote the corresponding quantizer:
$$G_k(x)=A_{k,i}\, {\rm if }\,x\in A_{k,i}.$$

The role of the feature vector in (\ref{eq1}) is now played by the
discrete quantity,
$$
B(k)=(G_k(X_{-\lambda _{k-1}}),\ldots,G_k(X_{-1}))=
G_k(X_{-\lambda_{k-1}}^{-1}).
$$
Now 
$$
{\tau}_k=
\min\{t>0 : G_k(X_{-\lambda_{k-1}-t}^{-1-t})=G_k(X_{-\lambda_{k-1}}^{-1})\}.
$$
Again, ergodicity implies that
$B(k)$ is almost surely to be found in the sequence
$G_k(X_{-\infty}^{-2})$, and
with this generalization of notation,
the $k$th estimate of $P(C|\bx ^-)$  is
still provided by formula (\ref{eq1}).

As in  Algoet's construct, the estimate $P_k$ is calculated from
observations of random size. Here the random sample size is $\lambda _k$.
To obtain a fixed sample size $t>0$ version,
let $\kappa_t$ be the maximum of integers $k$ for which $\lambda _k \le t$.
Put
\begin{equation}
\label{algdefpt}
\hat P_{-t}(C)=P_{\kappa_t}(C).
\end{equation}

\begin{theorem} \label{Theorem1} Under the stationary ergodic assumption
regarding $\{X_n\}$ and under the estimator
constructs (\ref{eq1}) and (\ref{algdefpt}) described above,
\begin{equation}
\label{eq3}
\lim_{k\rightarrow\infty} P_k(\cdot)= P(\cdot|\bx^-)\ \ {\rm a.s.,}
\end{equation}
and
\begin{equation}
\label{eq4}
\lim_{t\rightarrow \infty}\hat P_{-t}(\cdot)= P(\cdot|\bx^-)\ \ {\rm a.s.,}
\end{equation}
in the weak topology of distributions.
\end{theorem}

\proof
To begin with, assume that  for some $m$, $C\in \sigma({\p}_m)$.
The first chore is to show that  a.s.,
$$P_k(C)\rightarrow P(C|\bx ^-). $$
For $k>m$ we have that
\begin{eqnarray*}
\lefteqn{ P_k (C)-P(C|\bx ^-) }\\
&=& {1\over k}\sum_{1\le j\le m}
 [1_C(X_{-\tau _j})-P(X_{-\tau _j}\in C|G_{j-1}(X_{-\lambda_{j-1}}^{-1}))]\\
&+&{(k-m) \over k} {1\over (k-m)}\sum_{m<j\le k}
 [1_C(X_{-\tau _j})-P(X_{-\tau _j}\in C|G_{j-1}(X_{-\lambda_{j-1}}^{-1}))]\\
&+& {1\over k}\sum_{1\le j\le k}
 P(X_{-\tau_j}\in C|G_{j-1}(X_{-\lambda_{j-1}}^{-1}))
-P(C|\bx ^-)\\
&=&P1_k+{(k-m) \over k} P2_k+P3_k.
\end{eqnarray*}
Obviously,
$$P1_k\to 0\ \  {\rm a.s.}$$

\noindent
Toward mastering $P2_k$, one observes that
 $P2_k$ is an average of bounded martingale
differences. To see this note that
$\sigma (G_{j}(X_{-\lambda_{j}}^{-1}))$ $j=0,1,\dots$ is monotone
increasing, and that  $1_C(X_{-\tau _j})$   is measurable on
$\sigma (G_{j}(X_{-\lambda_{j}}^{-1}))$ for $j\ge m$. 
The convergence of $P2_k$  
can be established by L\'{e}vy's 
classical result, namely, the Ces\`{a}ro means of a bounded sequence 
of martingale differences converge to zero almost surely. 
For a version suited to our needs, 
see, for example, Theorem 3.3.1 in Stout (1974).    
One may even obtain rates  for $P2_k$ through the use
of Azuma's (1967) exponential bound for martingale differences. 
We have to prove that $$P3_k\to 0\ \ {\rm  a.s.}$$
By Lemma~\ref{lemma2} in the appendix,
$$
P(X_{-\tau _j}\in C|G_{j-1}(X_{-\lambda_{j-1}}^{-1}))= 
P(X_0\in C|G_{j-1}(X_{-\lambda_{j-1}}^{-1})).
$$
Using this we get
$$P3_k={1\over k}\sum_{1\le j\le k}
 P(X_{-\tau _j}\in C|G_{j-1}(X_{-\lambda_{j-1}}^{-1}))
-P(C|\bx ^-)$$
$$={1\over k}\sum_{1\le j\le k}
 P(X_0\in C|G_{j-1}(X_{-\lambda_{j-1}}^{-1}))
-P(C|\bx ^-).$$
By assumption,
$$\sigma (B(j))\uparrow \sigma(\bx ^-),$$
which implies that
$$\sigma(G_{j}(X_{-\lambda_{j}}^{-1}))\uparrow \sigma(\bx ^-).$$
Consequently by the a.s. martingale convergence theorem we have that
$$ P(X_0\in C|G_{j}(X_{-\lambda_{j}}^{-1}))\to P(C|\bx ^-)\ \
{\rm  a.s.,}$$
and thus by the Toeplitz lemma (cf. Ash (1972) )
$$P3_k\to 0\ \ {\rm  a.s.}$$

\noindent
Let D denote the countably infinite set of $x$'s for which $(-\infty ,x]\in
\sigma ({\p}_k)$ for sufficiently large $k$. By assumption, $D$ is dense in
$\R$.
Define
$$F_k(x)=P_k((-\infty ,x]).$$
Also, set
$$F(x)=P((-\infty ,x]|\bx ^-).$$
 By the preceding development we have the almost sure event $H$ such that
 on $H$ for all $x\in D$

\begin{equation}
\label{eq5b}
F_k(x)\to F(x).
\end{equation}
Since $D$ is dense in $\R$, we have (\ref{eq5b}) on $H$ and for all continuity points of
$F(\cdot)$, and (\ref{eq3}) is proved. 
The convergence (\ref{eq4}) is an obvious consequence of (\ref{eq3}).

\section{Estimation of auto-regression functions}
\label{sec3}

The next result uses estimators

\begin{equation}
\label{eqalgregk}
R_k={1\over k}\sum_{1\le j\le k} X_{-\tau _j}
\end{equation}
and
\begin{equation}
\label{eqalgregt}
\hat R_{-t}={1\over \kappa_t}\sum_{1\le j\le \kappa_t} X_{-\tau _j}.
\end{equation}
\begin{corollary} Assume that for some number $D$,
a.s., $|X_0|\le D<\infty.$
Under the stationary ergodic assumption
regarding $\{X_n\}$ and under the estimator
constructs (\ref{eqalgregk}) and (\ref{eqalgregt}) described above,
\begin{equation}
\label{eqregk}
\lim_{k\rightarrow\infty} R_k= E(X_0|\bx^-)\ \ {\rm a.s.,}
\end{equation}
and
\begin{equation}
\label{eqregt}
\lim_{t\rightarrow \infty} \hat R_{-t}= E(X_0|\bx^-) \ \ {\rm a.s.}
\end{equation}

\end{corollary}

\proof  Define the function
$$
\phi(x)= \left\{ \begin{array}{ll}
        D,\,  &\mbox{if $x>D$} \\
        x,\,  &\mbox{if $-D\le x \le D$}\\
        -D,\, &\mbox{if $x<-D$}
\end{array}
\right.
$$
Then
$$R_k=\int xP_k(dx)=\int\phi(x)P_k(dx)$$
$$\to \int\phi(x)P(dx|\bx^-)=\int xP(dx|\bx^-)=E(X_0|\bx^-). $$
because of Theorem~\ref{Theorem1} and the fact that convergence in distribution
implies the convergence of integrals of the bounded continuous function $\phi$
with respect to the actual distributions (Billingsley (1968)).
Thus the proof of (\ref{eqregk}) is complete. The proof
of (\ref{eqregt}) follows in the same way; just put $\hat P_{-t}$ in 
place of
$P_k$.

\bigskip

The estimates $\hat R_{-t}$ converge almost surely to $E(X_0|\bx^-)$ and are 
uniformly bounded so $|\hat R_{-t}-E(X_0|X^{-1}_{-t})|\to 0$ 
also in mean. Motivated by Bailey (1976), consider the estimator 
$\hat R_{t}(\omega)=\hat R_{-t}(T^t\omega)$ which is defined in terms of 
$(X_0,\ldots,X_{t-1})$ in the same way as $\hat R_{-t}(\omega)$ was defined 
in terms of $(X_{-t},\ldots,X_{-1})$. ($T$ denotes the left shift operator. )
The estimator $\hat R_t$ may be viewed as an on-line predictor of $X_t$. 
This predictor  has special significance not only because of potential
applications, but additionally because
Bailey (1976) proved that it is impossible to construct estimators $\hat R_t$ 
such that always $\hat R_t-E(X_t|X^{t-1}_0)\to 0$ almost surely. 
An immediate consequence of Corollary 1 is that convergence in
probability is verified.  
That is, the shift transformation $T$ is measure preserving hence convergence 
$\hat R_{-t}-E(X_0|X^{-1}_{-t})\to 0$ in $L^1$ implies convergence 
$\hat R_t-E(X_t|X^{t-1}_0)\to 0$ in $L^1$ and in probability.

\section{Pattern recognition}
\label{sec4}

Consider the $2$-class pattern recognition problem with $d$-dimensional feature
vector $X_0$ and binary valued label $Y_0$. Let
${\cal D}^-=(X_{-\infty}^{-1},Y_{-\infty}^{-1})$ be the data. In conventional pattern
recognition problems $(X_0,Y_0)$ and ${\cal D}^-$ are independent, so the best
possible decision based on $X_0$ and based on $(X_0,{\cal D}^-)$ are the same.
Here assume that $\{(X_i,Y_i)\}$ is a doubly infinite stationary and
ergodic sequence. The classification problem is to decide on $Y_0$ for
given data $(X_0,{\cal D}^-)$ in order to minimize the probability of misclassification.
The Bayes decision $g^*$ is the best possible one. Let $\eta(X_0,{\cal D}^-)$
be the {\it a posteriori} probability of $Y_0=1$ (regression function):
$$\eta (X_0,{\cal D}^-)= P(Y_0=1|X_0,{\cal D}^-)=E(Y_0|X_0,{\cal D}^-).$$
Then $g^*(X_0,{\cal D}^-)=1$ if $\eta (X_0,{\cal D}^-)\ge1/2$ and $0$
otherwise. For an arbitrary approximation $\eta_k=\eta_k(X_0,{\cal D}^-)$
put
$g_k=g_k(X_0,{\cal D}^-)=1$ if $\eta_k\ge1/2$ and $0$
otherwise. Then it is easy to see (cf. Devroye and Gy\"orfi (1985), Chapter 10)
that
\begin{eqnarray}
0&\le&
 P(g_k\ne Y_0|X_0,{\cal D}^-)
-P(g^*(X_0,{\cal D}^-)\ne Y_0|X_0,{\cal D}^-)\nonumber\\
&\le&\label{eq6}
2|\eta_k-\eta (X_0,{\cal D}^-)|.
\end{eqnarray}
The estimation is a slight modification of (\ref{eq1}).
Define the sequences $\lambda _{k-1}$ and $\tau _k$ recursively ($k=1,2,\dots$).
Put $\lambda _0=1$ and $\tau _k$ be the time between  the 
occurrence of the pattern
$$B(k)=(G_k(X_{-\lambda _{k-1}}),Y_{-\lambda _{k-1}},\ldots,
G_k(X_{-1}),Y_{-1},G_k(X_0))$$
at time $0$ and the last occurrence of the same pattern 
in ${\cal D}^-$. 
More precisely, 
$$
{\tau}_k=
\min\{t>0 : G_k(X_{-\lambda_{k-1}-t}^{-t})=G_k(X_{-\lambda_{k-1}}^{0}),
Y_{-\lambda_{k-1}-t}^{-1-t}=Y_{-\lambda_{k-1}}^{-1})\}.
$$
Put
$$\lambda _k=\tau _k+\lambda _{k-1}.$$
The observed vector
$B(k)$ a.s. takes a value of positive probability; thus by ergodicity
$B(k)$ has occurred with probability $1$.
One denotes the $k$th estimate of $\eta (X_0,{\cal D}^-)$ by
$\eta_k$,
and defines it to be
\begin{equation}
\label{algpatk}
\eta_k={1\over k}\sum_{1\le j\le k} Y_{-\tau _j}.
\end{equation}
\begin{corollary} Under the stationary ergodic assumption
regarding the process $\{(X_n,Y_n)\}$ and under the estimator
construct (\ref{algpatk}) described above,
\begin{equation}
\label{eq7}
P(g_k\ne Y_0|X_0,{\cal D}^-) \to
P(g^*(X_0,{\cal D}^-)\ne Y_0|X_0,{\cal D}^-) \ \ {\rm a.s.}
\end{equation}
\end{corollary}

\proof Because of (\ref{eq6}), we get (\ref{eq7}) from
$$\eta_k\to \eta (X_0,{\cal D}^-) \ \ {\rm a.s.,}$$
the proof of which is similar to the proof of Theorem~\ref{Theorem1}.

\noindent
{\bf Remark.}
It is also possible to construct a version of this estimate with
fixed sample size $t>0$ in the same way as in (\ref{algdefpt}) and
(\ref{eqalgregt}). 

\section{Appendix}
\label{sec6}

In the sequel, we use the notation of Section~\ref{sec2}.

\begin{lemma}
\label{lemma2}
Under the stationary ergodic assumption regarding
$\{X_n\}$, for $j=1,2,\dots$,
\[
P(X_{-\tau _j}\in C|G_{j-1}(X_{-\lambda_{j-1}}^{-1}))=
P(X_0\in C|G_{j-1}(X_{-\lambda_{j-1}}^{-1})).
\]
\end{lemma}
\proof
First of all, note that by definition,
\begin{eqnarray*}
\lefteqn{\sigma (G_{j-1}(X_{-\lambda_{j-1}}^{-1}))={\cal F}_{j-1}}\\
&=&\sigma (\{G_{j-1}(X_{-m}^{-1})=b_{-m}^{-1}, \lambda _{j-1}=m\}; 
b_{-m}^{-1}, m=1,2,\dots), 
\end{eqnarray*}
where $b_{-m}^{-1}$ is an $m$-vector of sets from the finite partition ${\p}_{j-1}$.

\noindent
Note also that
$$B=\{G_{j-1}(X_{-m}^{-1})=b_{-m}^{-1}, \lambda_{j-1}=m\}$$
are the (countable many) generating atoms  of
${\cal F}_{j-1}$, so
we have to show that for any atom $B$
the following equality holds:
$$ P(B\cap\{ X_{-\tau _j}\in C\}) 
= P(B\cap\{ X_0\in C\}).$$

\noindent
$\lambda_{j-1}$ is a stopping time, 
$B$ is an $m$-dimensional cylinder set, which means that 
$ b_{-m}^{-1}$ determines whether $\lambda_{j-1}\ne m$
(in which case $B=\emptyset$ and the statement is trivial) or 
$\lambda_{j-1}=m$ and then 
$$B=\{G_{j-1}(X_{-m}^{-1})=b_{-m}^{-1}\}.$$

\noindent 
For $j=1,2,\dots$ let 
$$
\tilde{\tau}_j=
\min\{0<t : G_{j}(X_{-\lambda_{j-1}+t}^{-1+t})=G_{j}(X_{-\lambda_{j-1}}^{-1})\}.
$$

\noindent
Now 
\begin{eqnarray*}
\lefteqn{T^{-l}[B\cap\{\tau_j=l,X_{-l}\in C\}]}\\
&=&T^{-l}[\{G_{j-1}(X_{-m}^{-1})=b_{-m}^{-1},G_{j}(X_{-m-l}^{-1-l})=G_{j}(X_{-m}^{-1}),\\
&& G_{j}(X_{-m-t}^{-1-t})\ne G_{j}(X_{-m}^{-1}), 0<t<l,X_{-l}\in C\}]\\
&=&\{G_{j-1}(X_{-m+l}^{-1+l})=b_{-m}^{-1},G_{j}(X_{-m}^{-1})=G_{j}(X_{-m+l}^{-1+l}),\\
&& G_{j}(X_{-m-t+l}^{-1-t+l})\ne G_{j}(X_{-m+l}^{-1+l}), 0<t<l,X_0\in C\}\\
&=&\{G_{j-1}(X_{-m+l}^{-1+l})=b_{-m}^{-1},G_{j}(X_{-m}^{-1})=G_{j}(X_{-m+l}^{-1+l}),\\
&& G_{j}(X_{-m+t}^{-1+t})\ne G_{j}(X_{-m+l}^{-1+l}), 0<t<l,X_0\in C\}\\
&=&\{G_{j-1}(X_{-m}^{-1})=b_{-m}^{-1},G_{j}(X_{-m}^{-1})=G_{j}(X_{-m+l}^{-1+l}),\\
&& G_{j}(X_{-m+t}^{-1+t})\ne G_{j}(X_{-m}^{-1}), 0<t<l,X_0\in C\}\\
&=&{B\cap\{\tilde\tau_j=l,X_0\in C\},}
\end{eqnarray*}
where $T$ denotes the left shift operator. 

\noindent
By stationarity, it follows that 
\begin{eqnarray*}
\lefteqn{P(B\cap \{X_{-\tau _j}\in C\})}\\
&=&\sum_{l=1}^{\infty}P(B\cap\{\tau_j=l,X_{-l}\in C\})\\
&=&\sum_{l=1}^{\infty}P(T^{-l}[B\cap\{\tau_j=l,X_{-l}\in C\}])\\
&=&\sum_{l=1}^{\infty}P(B\cap\{\tilde \tau_j=l,X_0\in C\})\\
&=&P(B\cap\{X_0\in C\}),
\end{eqnarray*}
and the proof of Lemma~\ref{lemma2} is complete.

\noindent
{\bf Acknowledgements }

The authors thank P. Algoet for 
his comments, suggestions, and encouragement.  Suggestions by the
referees have been helpful.
The second author's work has been supported, in part,
by NIH grant No. R01 A129426.

\end{document}